\documentclass[11pt]{amsart}
\usepackage{thm-restate,amsmath,amssymb,enumerate,graphicx}
\usepackage[
  colorlinks,
  linkcolor = blue,
  citecolor = blue,
  urlcolor = blue]{hyperref}

\newtheoremstyle{plainsl}%
	{\topsep}
	{\topsep}
	{\slshape} 
	{}
	{\normalfont\bfseries}
	{.}
	{ }
	{}
\swapnumbers

\newtheorem{theorem}{Theorem}[section]
\newtheorem{lemma}[theorem]{Lemma}
\newtheorem{corollary}[theorem]{Corollary}
\newtheorem{conjecture}[theorem]{Conjecture}
\newtheorem{definition}{Definition}[section]

\newcommand\cref[1]{Corollary~\ref{cor:#1}}

\renewcommand\proof{\noindent\textsl{Proof. }}
\newcommand\sqr[2]{{\vbox{\hrule height.#2pt
    \hbox{\vrule width.#2pt height#1pt \kern#1pt
        \vrule width.#2pt}\hrule height.#2pt}}}
\renewcommand\qed{%
	\ifmmode\eqno\sqr53
	\else\nolinebreak\ \hfill\sqr53\medbreak\fi}

\newcommand{\qarrow}{\xrightarrow{q}}
\newcommand{\bartheta}{\bar{\vartheta}}

\DeclareMathOperator{\tr}{Tr}

\DeclareMathOperator{\msum}{sum}
\newcommand{\cart}{\mathbin{\square}}

\begin{document}

\sloppy
\title[Sabidussi vs. Hedetniemi]{Sabidussi Versus Hedetniemi for Three Variations of the Chromatic Number}
\author[Godsil, Roberson, \v{S}\'{a}mal, and Severini]{Chris Godsil, David E.~Roberson, Robert \v{S}\'{a}mal, and Simone Severini}

\address{
        Chris Godsil, David Roberson, 
        Department of Combinatorics and Optimization, 
        University of Waterloo, 
        Waterloo, Canada}
\email{\{cgodsil,droberso\}@uwaterloo.math.ca}
\address{Robert \v{S}\'{a}mal, Computer Science Institute of Charles University, Prague, Czech Republic}
\email{samal@iuuk.mff.cuni.cz}
\address{
        Simone Severini\\
        Department of Computer Science\\
        University College London\\
        London, United Kingdom}
\email{s.severini@ucl.ac.uk}

\thanks{This research was supported by a grant from the
Natural Sciences and Engineering Research Council of Canada, the Royal Society, the Karel Jane\v{c}ek Science \& Research Endowment (NFKJ) grant 201201, grant LL1202
ERC CZ of the Czech Ministry of Education, Youth and Sports, and by grant GA \v{C}R P202-12-G061.}

\date{\today}
\maketitle

\begin{abstract}
We investigate vector chromatic number $(\chi_{vec})$, Lov\'{a}sz $\vartheta$-function of the complement $(\bartheta)$, and quantum chromatic number $(\chi_q)$ from the perspective of graph homomorphisms. We prove an analog of Sabidussi's theorem for each of these parameters, i.e.~that for each of the parameters, the value on the Cartesian product of graphs is equal to the maximum of the values on the factors. We also prove an analog of Hedetniemi's conjecture for $\bartheta$, i.e.~that its value on the categorical product of graphs is equal to the minimum of its values on the factors. We conjecture that the analogous results hold for vector and quantum chromatic number, and we prove that this is the case for some special classes of graphs.


\end{abstract}

\section{Introduction}

The chromatic number is a well known graph parameter which can be defined in terms of homomorphisms. A graph homomorphism from $G$ to $H$ is a function $\varphi: V(G) \to V(H)$ such that $\varphi(u)$ is adjacent to $\varphi(v)$ whenever $u$ is adjacent to $v$. In this terminology, a graph $G$ is $n$-colorable if and only if there exists a homomorphism from $G$ to $K_n$. There are many interesting variants of chromatic number which can also be defined via homomorphisms. We are concerned with the following:
\begin{itemize}
\item Vector chromatic number $(\chi_{vec})$
\item Strict vector chromatic number $(\bartheta)$
\item Quantum chromatic number $(\chi_q)$
\end{itemize}
As the notation suggests, the strict vector chromatic number is equal to the Lov\'{a}sz $\vartheta$-function of the complement, i.e.~$\bartheta(G) = \vartheta(\overline{G})$ where $\overline{G}$ denotes the complement of graph $G$. We, however, do not define it in this way, rather we approach both $\bartheta$ and $\chi_{vec}$ in terms of homomorphisms. From this viewpoint they can both be seen as relaxations of chromatic number defined in terms of assigning unit vectors to vertices such that vectors assigned to adjacent vertices have some specified inner product. Quantum chromatic number can also be viewed in terms of homomorphisms, however for this parameter we assign to vertices tuples of orthogonal projectors which must satisfy certain constraints.

In this paper we are concerned with how these parameters behave on certain graph products. We are in particular focused on the Cartesian and categorical products, denoted by $G \cart H$ and $G \times H$ respectively. A well known theorem of Sabidussi~\cite{Sabidussi} states that the chromatic number of the Cartesian product of two graphs is equal to the maximum of the chromatic numbers of its factors. An equally, if not more, well known conjecture of Hedetniemi proposes that the chromatic number of the categorical product of two graphs is equal to the minimum of the chromatic numbers of the factors. Our aim is to prove or make steps towards proving analogs of these two statements for the three parameters above. The rest of the paper is outlined as follows.



In Section~\ref{sec:basics} we define the basic concepts and notation used throughout the paper. We give the background needed for our results on vector and strict vector colorings in Section~\ref{subsec:vectcolor}. This is followed by Section~\ref{sec:bartheta}, in which we show that analogs of Sabidussi's theorem hold for $\bartheta$ and $\chi_{vec}$, and that the $\bartheta$ version of Hedetniemi's conjecture is true. In Section~\ref{sec:1homogeneous}, we investigate a class of graphs called 1-homogeneous graphs, which include edge transitive graphs. We give an explicit formula for $\bartheta$ and $\chi_{vec}$ for these graphs in terms of their largest and smallest eigenvalues. As a consequence, we see that these two parameters coincide for this class of graphs, and thus the $\chi_{vec}$ version of Hedetniemi's conjecture holds for 1-homogeneous graphs. In Section~\ref{subsec:quantum}, we introduce quantum homomorphisms and give the background needed for our results on quantum chromatic number. Then in Section~\ref{sec:quantum}, we prove the quantum analog of Sabidussi's theorem, and show that quantum Hedetniemi's conjecture holds for a family of graphs which initiated the study of quantum chromatic number. 


\section{Preliminaries}\label{sec:basics}

Here we give the background on the basic tools such as homomorphisms and graph products that we use throughout the paper. For a more detailed introduction we refer the reader to \cite{tardif, nesetril} for homomorphisms, and \cite{products} for graph products.

Let $G$ and $H$ be graphs (by which we mean undirected simple finite graphs). We denote the existence of a homomorphism from $G$ to $H$ by writing $G \rightarrow H$. It is easy to see that homomorphisms compose, so $G \rightarrow H \rightarrow K$ implies $G \rightarrow K$. In fact, graphs with homomorphisms form a category. More relevant for graph theory is that many graph theoretic notions can be simply expressed in terms of homomorphisms. In particular, a graph $G$ is $n$-colorable if and only if $G \rightarrow K_n$.

A graph parameter $f$ is called \emph{homomorphism-monotone} if we have $f(G) \le f(H)$ whenever $G \rightarrow H$. Examples of homomorphism-monotone parameters include $\chi, \chi_f, \omega$, etc. We will see that the three parameters, $\chi_{vec}$, $\bartheta$, and $\chi_q$ are homomorphism-monotone as well, and they are even \emph{quantum homomorphism-monotone}. (Quantum homomorphisms will be defined in Section~\ref{subsec:quantum}).

Given graphs $G$ and $H$, we define four graphs with vertex set $V(G) \times V(H)$. In the \emph{categorical product} $G \times H$ (also called direct, or tensor product), tuples $(u_1, v_1)$ and $(u_2, v_2)$ are adjacent if and only if
\[u_1 \sim u_2 \text{ and } v_1 \sim v_2.\]
In the \emph{Cartesian product} $G \cart H$, tuples $(u_1, v_1)$ and $(u_2, v_2)$ are adjacent if and only if
\[(u_1  \sim u_2 \text{ and } v_1 = v_2) \text{ or } (u_1 = u_2 \text{ and } v_1 \sim v_2).\]

The \emph{strong product} $G \boxtimes H$ is defined as the edge union of $G \times H$ and $G \cart H$. In the \emph{disjunctive product} $G * H$ (also referred to as the conormal product) the tuples $(u_1, v_1)$ and $(u_2, v_2)$ are adjacent if $u_1 \sim u_2$ or $v_1 \sim v_2$. It is trivial to see that $G \boxtimes H$ is a subgraph of $G * H$, and a little thought reveals that $\overline{G * H} = \overline{G} \boxtimes \overline{H}$, where $\overline{G}$ denotes the complement of $G$.

It is easy to see that $G \rightarrow G \cart H$ (and also $H \rightarrow G \cart H$), indeed, $G \cart H$ contains copies of both $G$ and $H$. By projecting onto each coordinate, we see that $G \times H \rightarrow G$ and $G \times H \rightarrow H$. (This is indeed true in any category and $G \times H$ is called the categorical product because it is, in fact, a product in the sense of category theory.) Consequently, we have the following lemma:

\begin{lemma}\label{lem:homomono}
If $f$ is a homomorphism-monotone graph parameter and $G$ and $H$ are graphs, then
\[f(G \times H) \le \min\{f(G), f(H)\} \ \text{ and } \ f(G \cart H) \ge \max\{f(G), f(H)\}.\]
\end{lemma}
\proof
This follows immediately from the fact that $G \times H \rightarrow G,H$ and $G,H \rightarrow G \cart H$.\qed

This lemma allows us to easily establish that $\chi_{vec}$, $\bartheta$, and $\chi_q$ must all satisfy the above inequalities. We can then ask if/when equality holds. Indeed, that is the main focus of this paper.

Much attention has specifically been given to the value of the chromatic number on the Cartesian and categorical products, and this is of course part of the motivation for our work. Applying the above lemma to the chromatic number, which is homomorphism-monotone, we obtain the following:
\begin{equation}\label{catprod}
\chi(G \times H) \le \min\{\chi(G), \chi(H)\}
\end{equation}
and
\begin{equation}\label{cartprod}
\chi(G \cart H) \ge \max\{\chi(G), \chi(H)\}.
\end{equation}
As mentioned above, a well known theorem of Sabidussi~\cite{Sabidussi} states that~(\ref{cartprod}) holds with equality, and we provide a proof here for the reader's convenience and comparison with Theorems~\ref{thm:sabibartheta} and~\ref{thm:qsabi}.

\begin{theorem}[Sabidussi 1957]\label{sabi}
For graphs $G$ and $H$,
\[\chi(G \cart H) = \max\{\chi(G), \chi(H)\}.\]
\end{theorem}
\proof
Let $m = \max\{\chi(G), \chi(H)\}$. Clearly, we need at least $m$ colors to color $G \cart H$. So it suffices to show that $G \cart H$ can be $m$-colored. There are colorings $g$ of $G$ and $h$ of $H$ using $m$ colors, which we may assume are the integers modulo $m$. It is easy to check that assigning $(g(u) + f(v)) \text{ mod } m$ to vertex $(u,v)$ gives an $m$-coloring of $G \cart H$.\qed

Determining whether (\ref{catprod}) holds with equality turns out to be much more difficult, and the following conjecture remains open to this day:
\begin{conjecture}[Hedetniemi 1966]\label{hedconj}
For all graphs $G$ and $H$,
\[\chi(G \times H) = \min\{\chi(G), \chi(H)\}.\]
\end{conjecture}

It is worth noting that many different versions of this statement have been either conjectured or proven since its inception. Perhaps most significantly, Zhu has recently proved in~\cite{frachedetniemi} that
\[\chi_f(G \times H) = \min\{\chi_f(G), \chi_f(H)\}\]
where $\chi_f$ denotes fractional chromatic number.

Zhu's proof makes use of the fact that fractional chromatic number can be written as a linear program and thus suffers from strong duality. As we will see below, the strong duality property of the semidefinite programs for $\bartheta$ and $\chi_{vec}$ are crucial for our proofs as well. This suggests that the lack of strong duality for chromatic number is one reason for the difficulty in attempting to prove Hedetniemi's conjecture.

\section{Vector and Strict Vector Colorings}\label{subsec:vectcolor}

Vector and strict vector colorings were first introduced in~\cite{chivec}, in which Karger, Motwani, and Sudan also defined the vector chromatic number. However, a parameter equal to $\chi_{vec}$ of the complement was actually introduced a few decades earlier in~\cite{MRR} and~\cite{Schrijver}, but this seems to have gone unnoticed by many. We will first define strict vector colorings:


\begin{definition}\label{def:strictvectcolor}
Let $\mathcal{S}^d$ denote the unit sphere in $\mathbb{R}^{d+1}$. For a graph $G$, a map $\varphi : V(G) \rightarrow \mathcal{S}^d$ is called a \emph{strict vector $k$-coloring} if whenever $u \sim v$,
\[\varphi(u)^T \varphi(v) = -\frac{1}{k-1}.\]
\end{definition}
So a strict vector $k$-coloring can be viewed as a homomorphism to the infinite graph whose vertices are unit vectors in $\mathbb{R}^d$ such that vectors $u$ and $v$ are adjacent whenever $u^Tv = -1/(k-1)$. The \emph{strict vector chromatic number} of $G$ is the infimum of real numbers $k$ such that $k > 1$ and $G$ admits a strict vector $k$-coloring (for all nonempty graphs this infimum can be obtained and thus is just the minimum). It has been shown~\cite{chivec} that the strict vector chromatic number of $G$ is equal to $\vartheta(\overline{G})$, and thus we use $\bartheta$ to denote this parameter. The Lov\'{a}sz $\vartheta$-function has been well studied and it possesses many interesting properties, some of which we will present below. For a more detailed look at this graph parameter we refer the reader to~\cite{Lovasz} and~\cite{knuth}.

If we relax the definition above to only require that adjacent vertices are assigned unit vectors which have inner product \emph{at most} $-1/(k-1)$, then we obtain what is known as a \emph{vector $k$-coloring}~\cite{chivec}. The smallest $k$ for which $G$ admits a vector $k$-coloring is the \emph{vector chromatic number} of $G$, and we denote this by $\chi_{vec}(G)$. The basic motivation behind these definitions is that mapping the vertices of the complete graph $K_n$ to the vertices of the $(n-1)$-dimensional simplex gives a (strict) vector $n$-coloring, and therefore any $n$-colorable graph is also (strict) vector $n$-colorable. This of course implies that $\chi_{vec}(G),\bartheta(G) \le \chi(G)$ for all graphs $G$. Though we will not, it can be proved (for instance by using the dual SDPs given below) that the above (strict) vector coloring of $K_n$ is optimal, and thus $\omega(G) \le \chi_{vec}(G),\bartheta(G)$.

Since any strict vector $k$-coloring is clearly a vector $k$-coloring as well, we have that $\chi_{vec}(G) \le \bartheta(G)$ for any graph $G$. Defining these parameters in terms of homomorphisms as above allows us to easily see that if both $\chi_{vec}$ and $\bartheta$ are homomorphism-monotone. Therefore, by Lemma~\ref{lem:homomono}, we have that
\begin{align*}
\chi_{vec}(G \times H) &\le \min\{\chi_{vec}(G), \chi_{vec}(H)\} \\
\bartheta(G \times H) &\le \min\{\bartheta(G), \bartheta(H)\}
\end{align*}
and
\begin{align*}
\chi_{vec}(G \cart H) &\ge \max\{\chi_{vec}(G), \chi_{vec}(H)\} \\
\bartheta(G \cart H) &\ge \max\{\bartheta(G), \bartheta(H)\}.
\end{align*}

It turns out that both $\bartheta$ and $\chi_{vec}$ can be written as semidefinite programs (SDPs). The practical advantage of this is that one can compute them to arbitrary precision in polynomial time. The other advantage is that we can use duality to assist in proving theorems. In general, strong duality does not hold for all SDPs, however one can show that it holds for the SDPs defining $\bartheta$ and $\chi_{vec}$ using Slater's condition. Below we give both the primal and dual SDPs for $\bartheta$ and $\chi_{vec}$. Here, $P \succeq 0$ means that the matrix $P$ is positive semidefinite, while $P \ge 0$ means that the entries of $P$ are nonnegative. We use $J$ to denote the all ones matrix, and $\circ$ to denote Schur product. The matrix $A$ in the SDP constraints refers to the adjacency matrix of $G$, and $\overline{A} := J - I - A$ is the adjacency matrix of $\overline{G}$.
\[
\begin{array}[t]{lcc}
 & \text{PRIMAL} & \text{DUAL} \vspace{.05in} \\
\bartheta(G) & \begin{array}[t]{ll}
\min & \lambda \\
\text{s.t.} & M \circ I = (\lambda - 1)I \\
 & M \circ A = -A \\
 & M \succeq 0
\end{array} & \begin{array}[t]{ll}
\max & \tr(JP) \\
\text{s.t.} & P \circ \overline{A} = 0 \\
 & \tr(P) = 1 \\
 & P \succeq 0
\end{array} \\
 & & \\
\chi_{vec}(G) & \begin{array}[t]{ll}
\min & \lambda \\
\text{s.t.} & M \circ I = (\lambda - 1)I \\
 & M \circ A \le -A \\
 & M \succeq 0
\end{array} & \begin{array}[t]{ll}
\max & \tr(JP) \\
\text{s.t.} & P \circ \overline{A} = 0 \\
 & \tr(P) = 1 \\
 & P \ge 0 \\
 & P \succeq 0
\end{array}
\end{array}
\]
To see that the SDPs for $\chi_{vec}$ are equivalent to the vector coloring definition of this parameter, one one can use the fact that the positive semidefinite matrix $M$ in the primal SDP is a Gram matrix of a set of vectors. Assigning these (normalized) vectors to the vertices of the graph gives a valid vector coloring of the appropriate value. The reverse procedure converts a vector coloring to a feasible solution to the primal. The same technique works for $\bartheta$ as well~\cite{Lovasz, knuth, Schrijver}.

%
%
%
%

\section{Strict Vector Chromatic Number}\label{sec:bartheta}

To prove the strict vector chromatic number version of Sabidussi's theorem, we need the following lemma which shows that any graph $G$ which can be strict vector $k$-colored, can also be strict vector $k'$-colored for any $k' \ge k$. For chromatic number, as well as vector chromatic number, this is trivial, since any $k$-coloring can be viewed as a $k'$-coloring for any $k' \ge k$.

\begin{lemma}\label{biggerk}
Suppose $G$ is a graph such that $\bartheta(G) = k$. Then for every real $k' \ge k$, there is a strict vector $k'$-coloring of $G$.
\end{lemma}
\proof
Let $\varphi : V(G) \rightarrow \mathbb{R}^d$ be a strict vector $k$-coloring of $G$. Let $t = -1/(k-1)$, and $t' = -1/(k'-1)$. As $k' \ge k > 1$, we have that $t \le t' < 0$. Consequently, there exists an $\alpha \in [0,1]$ such that $\alpha^2 t + (1 - \alpha^2) = t'$. Define the mapping $\varphi' = (\alpha \varphi, \sqrt{1 - \alpha^2})$. It is easy to check that $\varphi'$ is a strict vector $k'$-coloring of $G$.\qed

We are now able to prove that Sabidussi's theorem holds for $\bartheta$.

\begin{theorem}\label{thm:sabibartheta}
For graphs $G$ and $H$,
\[\bartheta(G \cart H) = \max\{\bartheta(G), \bartheta(H)\}.\]
\end{theorem}
\proof
As we have already seen in Section~\ref{subsec:vectcolor},
\[\bartheta(G \cart H) \ge \max\{\bartheta(G), \bartheta(H)\}.\]
Thus, we only need to show the reverse inequality. Let $k = \max\{\bartheta(G), \bartheta(H)\}$. By Lemma~\ref{biggerk}, there exist strict vector $k$-colorings $g: V(G) \to \mathbb{R}^{d_1}$ and $h: V(H) \to \mathbb{R}^{d_2}$. We will consider the tensor product $g \otimes h : V(G \cart H) \to \mathbb{R}^{d_1 d_2}$. Explicitly, we put $(g \otimes h)(u,v) = g(u) \otimes h(v)$, where $u \in V(G)$ and $v \in V(H)$.

Now consider an edge of the form $(u,v) (u',v)$ in $G \cart H$. Let $t = -1/(k-1)$ as in Definition~\ref{def:strictvectcolor}. Using standard properties of the tensor product we get 
\[\left(g(u) \otimes h(v) \right)^T \left(g(u') \otimes h(v) \right) = \left(g(u)^Tg(u')\right) \left(h(v)^T h(v)\right) = t \cdot 1 = t.\]
By symmetry, we get the same condition for edges of the form $(u,v)(u,v')$. Consequently, $g \otimes h$ is a strict vector $k$-coloring of $G \cart H$, as required.\qed

We also have the following:

\begin{lemma}
For graphs $G$ and $H$,
\[\chi_{vec}(G \cart H) = \max\{\chi_{vec}(G), \chi_{vec}(H)\}.\]
\end{lemma}
\proof
Same as in Theorem~\ref{thm:sabibartheta}, without the need for Lemma~\ref{biggerk} since vector $k$-colorings only require that adjacent vertices have inner product \emph{at most} $-1/(k-1)$.\qed

We will use Lemma~\ref{thm:sabibartheta} to prove the $\bartheta$ version of Hedetniemi's conjecture, but we will also need some basic facts about how $\bartheta$ behaves on the strong and disjunctive products, as well as the edge union of two graphs.

In~\cite{Lovasz} it was shown that $\vartheta(G \boxtimes H) = \vartheta(G)\vartheta(H)$. A slight modification of the same proof shows that $\vartheta(G * H) = \vartheta(G)\vartheta(H)$, and in fact this is proven in~\cite{knuth}. Translating these two facts into terms of $\bartheta$, we obtain the following lemma.

\begin{lemma}\label{lem:barthetastrongprod}
For graphs $G$ and $H$,
\[\bartheta(G \boxtimes H) = \bartheta(G)\bartheta(H) = \bartheta(G * H).\]
\end{lemma}
\proof
Since $\overline{G * H} = \overline{G} \boxtimes \overline{H}$ (and equivalently $\overline{G \boxtimes H} = \overline{G} * \overline{H}$), we have that
\[\bartheta(G * H) = \vartheta(\overline{G} \boxtimes \overline{H}) = \vartheta(\overline{G})\vartheta(\overline{H}) = \bartheta(G)\bartheta(H)\]
and
\[\bartheta(G \boxtimes H) = \vartheta(\overline{G} * \overline{H}) = \vartheta(\overline{G})\vartheta(\overline{H}) = \bartheta(G)\bartheta(H).\]\qed

From this lemma we can easily obtain the following corollary which is analogous to a well known upper bound on the chromatic number of the union of two graphs. Given two graphs $G$ and $H$ on the same vertex set $V$, we use $G \cup H$ to denote the graph with vertex set $V$ and edge set $E(G) \cup E(H)$.

\begin{corollary}\label{cor:barthetaunion}
If $G$ and $H$ are graphs on the same vertex set $V$, then
\[\bartheta(G \cup H) \le \bartheta(G)\bartheta(H).\]
\end{corollary}
\proof
The vertices of the form $(v,v)$ for $v \in V$, induce a subgraph of $G * H$ isomorphic to $G \cup H$ and thus
\[\bartheta(G \cup H) \le \bartheta(G * H) = \bartheta(G)\bartheta(H)\]
by Lemma~\ref{lem:barthetastrongprod}.\qed

With these tools in hand, we are able to give a simple and elegant proof of Hedetniemi's conjecture for $\bartheta$.

\begin{theorem}\label{thm:hedbartheta}
For graphs $G$ and $H$,
\[\bartheta(G \times H) = \min\{\bartheta(G), \bartheta(H)\}.\]
\end{theorem}
\proof
We have already seen that
\[\bartheta(G \times H) \le \min\{\bartheta(G), \bartheta(H)\}.\]
So we only need to show the reverse inequality. For this we observe that $G \boxtimes H = (G \cart H) \cup (G \times H)$. Using Corollary~\ref{cor:barthetaunion} for $G \cart H$ and $G \times H$, as well as Lemma~\ref{lem:barthetastrongprod}, we obtain
\[\bartheta(G) \bartheta(H) = \bartheta(G \boxtimes H) \le \bartheta(G \cart H) \bartheta(G \times H).\]
Combining this with Theorem~\ref{thm:sabibartheta} finishes the proof.\qed

\section{1-Homogeneous Graphs}\label{sec:1homogeneous}

A graph $G$ is said to be \emph{1-homogeneous} if it satisfies the following two conditions:
\begin{enumerate}
\item The number of closed walks of length $k$ in $G$ that begin at a vertex $u$ is independent of $u$ for all $k \in \mathbb{Z}$.

\item The number of walks of length $k$ in $G$ that begin at vertex $u$ and end at adjacent vertex $v$ is independent of the edge $uv$.
\end{enumerate}

The first condition can be viewed as a type of combinatorial relaxation of vertex transitivity. Indeed, it is easy to see that any vertex transitive graph has this property. The second condition can similarly be viewed as a combinatorial relaxation of edge transitivity, and again, any edge transitive graph trivially has this property. So any graph which is both edge and vertex transitive is 1-homogeneous. Note that letting $k = 2$ in the first condition guarantees any such graph is regular.

Though 1-homogeneous graphs are not a well known class of graphs, they include several well known classes of graphs. In particular, distance regular (and thus strongly regular) graphs are 1-homogeneous. ore generally, any graph which is a single class in an association scheme is 1-homogeneous. We will also see that any edge transitive graph is either 1-homogeneous or bipartite, thus the results of this section apply to all of these classes of graphs.

If $A$ is the adjacency matrix of a graph $G$, then the $uv$-entry of $A^k$ is the number of walks of length $k$ in $G$ starting at $u$ and ending at $v$. From this it is easy to see that $G$ being 1-homogeneous is equivalent to the existence of constants $b_k$ and $c_k$ for all $k \in \mathbb{N}$ such that
\begin{equation}\label{eqn:1hom}
A^k \circ I = b_k I \ \ \ \& \ \ \ A^k \circ A = c_kA.
\end{equation}

In this section, we will present an explicit formula for the vector chromatic number of a 1-homogeneous graph in terms of its largest and smallest eigenvalues. Furthermore, we will show that $\bartheta$ and $\chi_{vec}$ are equal in this case. As a result, we will see that the vector chromatic number version of Hedetniemi's conjecture holds for all 1-homogeneous graphs. The results of this section rely heavily on the SDP formulations of $\chi_{vec}$ and $\bartheta$ given in Section~\ref{subsec:vectcolor}.

Before we give our results on 1-homogeneous graphs, we prove a general lower bound on vector chromatic number. The following two lemmas are from a set of unpublished notes of Godsil~\cite{interestinggraphs}.

\begin{lemma}
Let $G$ be a graph with $n$ vertices, $e$ edges, and least eigenvalue~$\tau$. Then
\[\chi_{vec}(G) \ge 1 - \frac{2e/n}{\tau}.\]
\end{lemma}
\proof
Let $A$ be the adjacency matrix of $G$. Then $A - \tau I \succeq 0$ and
\[(A - \tau I) \circ \overline{A} = 0.\]
Since $\tau < 0$, we have that
\[A - \tau I \ge 0.\]
Furthermore, since
\[\tr(A - \tau I) = -n\tau,\]
the matrix $\frac{1}{-n\tau}(A - \tau I)$ is a feasible solution to the dual formulation of $\chi_{vec}$ with objective value $1 - \frac{2e/n}{\tau}$. This gives the lower bound and proves the lemma.\qed

Note that $2e/n$ is the average degree of the graph, and is thus simply the degree for regular graphs. The next lemma states that the above bound is tight for 1-homogeneous graphs.

\begin{lemma}\label{1chivec}
If $G$ is 1-homogeneous with degree $k$ and least eigenvalue $\tau$, then
\[\chi_{vec}(G) = \bartheta(G) = 1 - \frac{k}{\tau}.\]
\end{lemma}
\proof
We make use of the identity $\tr(A^T B) = \msum(A \circ B)$ where $\msum(M)$ denotes the sum of all the entries of the matrix $M$. From the previous lemma, we have that $\chi_{vec}(G) \ge 1 - \frac{2e/n}{\tau} = 1 - \frac{k}{\tau}$. As we saw in Section~\ref{subsec:vectcolor}, $\chi_{vec}(G) \le \bartheta(G)$, and thus we only need to show that $\bartheta(G) \le 1 - \frac{k}{\tau}$. To do this we will find a suitable solution to the primal SDP formulation of $\bartheta$. Let $A$ be the adjacency matrix of $G$ and let $E_\tau$ denote the projection onto the $\tau$-eigenspace of $A$. Since $G$ is 1-homogeneous and $E_\tau$ is a polynomial in $A$, by~(\ref{eqn:1hom}) we have that
\[E_\tau \circ I = bI \ \text{ and } \ E_\tau \circ A = cA,\]
for some constants $b$ and $c$. Let $r$ be the rank of $E_\tau$. Since $E_\tau$ is a projection, $\tr(E_\tau) = r$, and so $b = r/n$. Now
\[\msum(A \circ E_\tau) = \tr(AE_\tau) = \tr(\tau E_\tau) = r\tau\]
and also
\[\msum(A \circ E_\tau) = c \cdot \msum(A) = cnk,\]
whence $c = r\tau/nk$. If we define
\[M := -\frac{nk}{r\tau}E_\tau,\]
then
\[M \circ I = -\frac{nk}{r\tau}bI = -\frac{k}{\tau}I = \left(\left(1 - k/\tau\right) - 1\right)I\]
and
\[M \circ A = -\frac{nk}{r\tau}cA = -A.\]
Since $E_\tau$ is a projection, $M$ is positive semidefinite, and is thus a feasible solution to the primal SDP with objective value $1 - \frac{k}{\tau}$. Therefore $\bartheta(X) \le 1 - \frac{k}{\tau}$ and the lemma is proven.\qed

Applying the above to edge transitive graphs, we obtain the following corollary
\begin{corollary}\label{cor:edgetheta}
If $G$ is edge transitive with greatest and least eigenvalues $\lambda$ and $\tau$ respectively, then
\[\chi_{vec}(G) = \bartheta(G) = 1 - \frac{\lambda}{\tau}.\]
\end{corollary}
\proof
First note that we can assume that $G$ has no isolated vertices since removing them does not change any of $\lambda$, $\tau$, $\chi_{vec}(G)$, or $\bartheta(G)$. If $G$ is also vertex transitive, then it is 1-homogeneous and $\lambda$ is the degree, $k$, of $G$. Thus the result holds by the above. If $G$ is not vertex transitive, then by Lemma~3.2.1 from~\cite{AGT}, it is nonempty and bipartite and thus $\chi_{vec}(G) = 2 = \bartheta(G)$. However, the spectrum of any bipartite graph is symmetric about zero, and so $\tau = -\lambda$, and therefore $1 - \frac{\lambda}{\tau} = 2$.\qed

Lemma~\ref{1chivec} above is the key to proving that the $\chi_{vec}$ analog of Hedetniemi's conjecture holds for 1-homogeneous graphs. But to be able to apply it we need the following lemma from~\cite{equiarboreal}:

\begin{lemma}[Godsil]\label{lem:1times}
If $G$ and $H$ are 1-homogeneous graphs, then the graph $G \times H$ is 1-homogeneous.\qed
\end{lemma}

We are now ready to prove the $\chi_{vec}$ version of Hedetniemi's conjecture for 1-homogeneous graphs.

\begin{theorem}\label{thm:1hedchivec}
If $G$ and $H$ are 1-homogeneous, then
\[\chi_{vec}(G \times H) = \min\{\chi_{vec}(G), \chi_{vec}(H)\}.\]
\end{theorem}
\proof
Since $G$ and $H$ are 1-homogeneous, and by Lemma~\ref{lem:1times} the product $G \times H$ is as well, we have that
\[\chi_{vec}(G \times H) = \bartheta(G \times H) = \min\{\bartheta(G), \bartheta(H)\} = \min\{\chi_{vec}(G), \chi_{vec}(H)\},\]
where the second equality follows from Theorem~\ref{thm:hedbartheta}.\qed

Note that one can also prove the above theorem without the aid of Theorem~\ref{thm:hedbartheta} by writing the degree and smallest eigenvalue of $G \times H$ in terms of the same parameters for $G$ and $H$.

As a corollary, we get that Hedetniemi's conjecture for $\chi_{vec}$ also holds for all edge transitive graphs.

\begin{corollary}\label{cor:hededge}
If $G$ and $H$ are edge transitive, then
\[\chi_{vec}(G \times H) = \min\{\chi_{vec}(G), \chi_{vec}(H)\}.\]
\end{corollary}
\proof
As in Corollary~\ref{cor:edgetheta}, we can assume that neither $G$ nor $H$ contains any isolated vertices. If both graphs are vertex transitive, then they are 1-homogeneous and the result holds by Theorem~\ref{thm:1hedchivec}. Otherwise, at least one of them is bipartite and therefore their product is bipartite and the conjecture holds trivially in this case.\qed

\section{Quantum Colorings}\label{subsec:quantum}

As a result of the continuing attempt to isolate the differences between quantum and classical mechanics, a large literature has developed which is devoted to the study of communication protocols based on the use of quantum resources, such as shared physical systems. In order to approach this problem with quantitative techniques and from a combinatorial angle, quantum colorings and the quantum chromatic number were introduced in \cite{Galliard02, Cleve04} and \cite{Avis} respectively. These concepts were further investigated in~\cite{qchrom, qchrom-aqis, SS12, MSS}.

A seminal result in the study of quantum colorings was the discovery of a family of graphs $\{\Omega_{4n} : n \in \mathbb{N}\}$ which exhibit an exponential separation between $\chi(\Omega_{4n})$ and $\chi_q(\Omega_{4n})$~\cite{BCW98, BCT99, Galliard02}. Here, $\Omega_n$ is the graph with vertex set $\{\pm 1\}^n$ such that orthogonal vectors are adjacent. In~\cite{frankrod} it was shown that when $n$ is a multiple of four, the graph $\Omega_n$ has chromatic number exponential in $n$. In contrast, it was shown in~\cite{BCW98, BCT99, Galliard02} that when $n$ is a power of two, $\chi_q(\Omega_n) \le n$. This result was extended to all $n$ divisible by four in~\cite{Avis}.

In~\cite{qhomos}, Man\v{c}inska and Roberson introduce the notion of quantum homomorphisms, which generalize quantum colorings in the same way that homomorphisms generalize colorings. It is this framework we will use to study quantum colorings and quantum chromatic number.
For a more detailed look at quantum colorings and quantum homomorphisms we refer the reader to~\cite{qchrom} and~\cite{qhomos}.

Though quantum homomorphisms were originally defined via a game played between two players and a referee, by the results of~\cite{qchrom, qhomos}, one can equivalently define them using homomorphisms. To do this, we require the following definition which comes from~\cite{qhomos}:

\begin{definition}
For a graph $G$ and integer $d$, let $M(G,d)$ be the the following graph. The vertices of $M(G,d)$ are the tuples ${\bf{E}} = (E_v)_{v \in V(G)}$ such that $E_v \in \mathbb{C}^{d \times d}$ is an orthogonal projector for all $v \in V(G)$ and
\begin{equation}\label{eq:meas}
\sum_{v \in V(G)} E_v = I.
\end{equation}
Two vertices ${\bf{E}} = (E_v)_{v \in V(G)}$ and ${\bf{E'}} = (E'_v)_{v \in V(G)}$ are adjacent if whenever $v \not\sim v'$,
\[E_v E'_{v'} = 0.\]
We refer to the graph $M(G,d)$ as the \emph{measurement graph of $G$ in dimension~$d$}.
\end{definition}

Note that in the above we do not consider a vertex to be adjacent to itself and thus $v \not\sim v$ for all $v \in V(G)$. Furthermore, note that condition~(\ref{eq:meas}) implies that for distinct vertices $v,v' \in V(G)$, we have that $E_vE_{v'} = 0$.

The reasoning behind the name of the measurement graph is that its vertices are what are known as ``projective quantum measurements". In general, a quantum measurement can consist of any positive semidefinite operators which sum to identity, but if each of the operators is a projection, then it is referred to as a projective measurement.

We say that $G$ has a quantum homomorphism to $H$, and write $G \qarrow H$, if $G \to M(H,d)$ for some $d \in \mathbb{N}$. We will also refer to a homomorphism from $G$ to $M(H,d)$ as a quantum homomorphism from $G$ to $H$. Note that if $\varphi$ is a homomorphism from $G$ to $H$, then the map which takes $u \in V(G)$ to the tuple whose $\varphi(u)$ coordinate is $I$ and all other coordinates are 0 is a quantum homomorphism from $G$ to $H$. Therefore $G \to H \Rightarrow G \qarrow H$.

Now that we have defined quantum homomorphisms, we can define quantum colorings and quantum chromatic number in the obvious way: a \emph{quantum $n$-coloring} of a graph $G$ is simply a quantum homomorphism from $G$ to $K_n$, and the \emph{quantum chromatic number} of $G$, denoted $\chi_q(G)$, is the minimum $n$ such that $G \qarrow K_n$. Note that since $G \to H \Rightarrow G \qarrow H$, for all $G$ and $H$, we have that $\chi_q(G) \le \chi(G)$ for all graphs $G$.

The definition of quantum homomorphism may seem a bit arbitrary, but it arises from the following physical considerations.

For graphs $G$ and $H$, the $(G,H)$-homomorphism game consists of two players, Alice and Bob, trying to convince a referee that they have a homomorphism from $G$ to $H$. More precisely, the referee sends Alice and Bob vertices $u_A, u_B \in V(G)$ respectively, and they respond with vertices $v_A, v_B \in V(H)$ accordingly. To win, the following conditions must be satisfied:
\begin{align*}
\text{if } u_A = u_B, \text{ then } v_A = v_B; \\
\text{if } u_A \sim u_B, \text{ then } v_A \sim v_B.
\end{align*}
Players can decide upon a strategy beforehand, but cannot communicate once play has commenced. The game is played for only one round, but we require a ``winning" strategy
to win with probability 1. It is not too difficult to see that classical players (who can use probabilistic strategies and have access to shared randomness) can win the
$(G,H)$-homomorphism game with certainty if and only if there exists a homomorphism from $G$ to~$H$. However, if players are allowed to perform quantum measurements on a
shared entangled state, then it is sometimes possible for them to win the $(G,H)$-homomorphism game even when $G \not\rightarrow H$. A general introduction to the theory of quantum entanglement can be found in~\cite{nielsen}. In Chapter~10 of~\cite{arorabarak} the interested reader may find a short elementary analysis of a different communication game, which exhibits an analogous difference between the classical and the quantum version. 

In~\cite{qchrom} it was proven that for $H = K_n$, the $(G,H)$-homomorphism game can be won by quantum players if
and only if $G \to M(H,d)$ for some $d \in \mathbb{N}$ (though it was not phrased in this way). In~\cite{qhomos} they note that the same proof works for any graph $H$ and
they introduce the measurement graph.

The general idea behind the correspondence between winning quantum strategies for the $(G,H)$-homomorphism game and homomorphisms from $G$ to $M(H,d)$ is as follows: Let $\varphi$ be a homomorphism from $G$ to $M(H,d)$. If Alice and Bob receive $u_A, u_B \in V(G)$ respectively, then Alice and Bob can perform measurements $\varphi(u_A)$ and $\varphi(u_B)^T$ on what is known as a ``maximally entangled state" to win the game. Here, $\varphi(u_B)^T$ corresponds to taking the transpose of each coordinate of $\varphi(u_B)$. The adjacency condition for $M(H,d)$ will correspond to the probability of outputting an incorrect response being zero.

%

In many ways quantum homomorphisms behave similarly to homomorphisms. In~\cite{qhomos} it was shown that they are transitive, i.e.~if $G \qarrow H \qarrow K$, then $G \qarrow K$. This means that $\chi_q$ is quantum homomorphism-monotone, i.e.~that $G \qarrow H \Rightarrow \chi_q(G) \le \chi_q(H)$. Note that in general $G \qarrow H$ does not imply that $\chi(G) \le \chi(H)$. Similarly, many other graph parameters defined via homomorphisms are not quantum homomorphism-monotone. However, in~\cite{qhomos} it was shown that both $\chi_{vec}$ and $\bartheta$ are quantum homomorphism-monotone, i.e.~$G \qarrow H$ implies that
\[\chi_{vec}(G) \le \chi_{vec}(H) \text{ and } \bartheta(G) \le \bartheta(H).\]
Since $\bartheta(K_n) = n$, it follows that strict vector chromatic number lower bounds quantum chromatic number. From this we obtain the following lemma:
\begin{lemma}\label{lem:chain}
For any graph $G$, we have
\[\chi_{vec}(G) \le \bartheta(G) \le \chi_q(G).\qed\]
\end{lemma}

We mentioned above that $\chi_q$ is quantum homomorphism-monotone. This is in fact a stronger condition than being homomorphism-monotone. Indeed, $G \rightarrow H \Rightarrow G \qarrow H \Rightarrow f(G) \le f(H)$ for any quantum homomorphism-monotone parameter $f$. Therefore $\chi_q$ is homomorphism-monotone, and thus by Lemma~\ref{lem:homomono} we have
\[\chi_q(G \cart H) \ge \max\{\chi_q(G), \chi_q(H)\},\]
and
\[\chi_q(G\times H) \le \min\{\chi_q(G), \chi_q(H)\}.\]
So we have seen that the easy directions of Sabidussi's theorem and Hedetniemi's conjecture hold for all three of the parameters we are investigating.

\section{Quantum Chromatic Number}\label{sec:quantum}

Here we will prove the quantum analog of Sabidussi's theorem, and use Theorem~\ref{thm:hedbartheta} to show that the quantum analog of Hedetniemi's conjecture holds in certain cases. First, we need the following lemma. We denote by $G[H]$ the lexicographic product of $G$ with $H$, for a definition see~\cite{products}.

\begin{lemma}\label{lem:qhomocart}
Suppose that $G,H,F,K$ are graphs such that $G \qarrow F$ and $H \qarrow K$. Then the following hold
\begin{enumerate}
\item $G \cart H \qarrow F \cart K$;
\item $G \times H \qarrow F \times K$;
\item $G \boxtimes H \qarrow F \boxtimes K$;
\item $G * H \qarrow F * K$;
\item $G[H] \qarrow F[K]$.
\end{enumerate}
\end{lemma}
\proof
We only give the proof for item (1), but it is obvious that a similar proof works for the others. For a function $f$ from vertices to tuples, we will use $f_u(v)$ to denote the $u$ coordinate of $f(v)$. Suppose that $\varphi^1$ and $\varphi^2$ are homomorphisms from $G$ to $M(F,d_1)$ and from $H$ to $M(K,d_2)$ respectively. Define $\varphi : V(G \cart H) \to V(M(F \cart K,d_1d_2))$ as follows:
\[\varphi_{(w,z)}(u,v) = \varphi^1_w(u) \otimes \varphi^2_z(v)\]
for all $(u,v) \in V(G \cart H)$ and $(w,z) \in V(F \cart K)$. First, we must show that $\varphi$ is indeed a map to the vertices of $M(F \cart K,d_1d_2)$. Since $\varphi^1_w(u)$ and $\varphi^2_z(v)$ are orthogonal projectors in dimensions $d_1$ and $d_2$ respectively, their tensor product is an orthogonal projector in dimension $d_1d_2$. Furthermore, since
\[\sum_{w \in V(F)} \varphi^1_w(u) = I \text{ for all } u \in V(G)\]
and
\[\sum_{z \in V(K)} \varphi^2_z(v) = I \text{ for all } v \in V(H)\]
we have that
\begin{align*}
\sum_{(w,z) \in V(F \cart K)} \varphi_{(w,z)}(u,v) &= \sum_{w \in V(F), z \in V(K)} \varphi^1_w(u) \otimes \varphi^2_z(v)  \\
&= \left(\sum_{w \in V(F)} \varphi^1_w(u)\right) \otimes \left(\sum_{z \in V(K)} \varphi^2_z(v)\right) \\
&= I \otimes I = I.
\end{align*}
Now recall from the definition of $M(F,d_1)$ that for $u \sim u' \in V(G)$, we have that $\varphi^1_w(u) \varphi^1_{w'}(u') = 0$ whenever $w \not\sim w'$. We also have that $\varphi^1_w(u) \varphi^1_{w'}(u) = 0$ for distinct $w,w' \in V(F)$, and the analogous conditions for $\varphi^2$. 

To show that $\varphi$ is a homomorphism, we must show that for $(u,v) \sim (u',v')$, we have $\varphi_{(w,z)}(u,v) \varphi_{(w',z')}(u',v') = 0$ whenever $(w,z) \not\sim (w'z')$. Since 
\[\varphi_{(w,z)}(u,v) \varphi_{(w',z')}(u',v') = \varphi^1_w(u) \varphi^1_{w'}(u') \otimes \varphi^2_z(v) \varphi^2_{z'}(v'),\]
it suffices to show that either $\varphi^1_w(u) \varphi^1_{w'}(u') = 0$ or $\varphi^2_z(v) \varphi^2_{z'}(v') = 0$.

Since $(u,v) \sim (u',v')$, without loss of generality we have that $u \sim u'$ and $v = v'$. The latter implies that $\varphi^2_z(v) \varphi^2_{z'}(v') = 0$ unless $z = z'$. However, if $z = z'$ and $(w,z) \not\sim (w',z')$, then we must have that $w \not\sim w'$ and thus $\varphi^1_w(u) \varphi^1_{w'}(u') = 0$. Therefore, we have shown that $\varphi$ is a homomorphism from $G \cart H$ to $M(F \cart K, d_1d_2)$, and thus $G \cart H \qarrow F \cart K$.\qed

%
%

We will in fact only need item (1) from the above lemma. We state the others simply because they follow from an essentially identical proof. Recall from Section~\ref{subsec:quantum} that quantum homomorphisms are transitive, and that $G \rightarrow H \Rightarrow G \qarrow H$ for any graphs $G$ and $H$. With these facts and the above lemma, we can easily prove the quantum version of Sabidussi's theorem.

\begin{theorem}\label{thm:qsabi}
For graphs $G$ and $H$,
\[\chi_q(G \cart H) = \max\{\chi_q(G), \chi_q(H)\}.\]
\end{theorem}
\proof
We saw in Section~\ref{subsec:quantum} that $\chi_q(G \cart H) \ge \max\{\chi_q(G), \chi_q(H)\}$, so we only need to show the other inequality. Let $n = \max\{\chi_q(G), \chi_q(H)\}$. Then we have that $G \qarrow K_n$ and $H \qarrow K_n$. Therefore, by Lemma~\ref{lem:qhomocart} and the original Sabidussi's theorem, we have
\[G \cart H \qarrow K_n \cart K_n \rightarrow K_n\]
and thus
\[G \cart H \qarrow K_n.\]
Therefore $\chi_q(G \cart H) \le n$, and we are done.\qed

Although we are not able to prove the general quantum version of Hedetniemi's conjecture, we can use the $\bartheta$ version of Hedetniemi's conjecture to prove a special case.

\begin{theorem}
Suppose that graphs $G$ and $H$ are such that $\chi_q(G) = \bartheta(G)$ and $\chi_q(H) = \bartheta(H)$. Then
\[\chi_q(G \times H) = \min\{\chi_q(G), \chi_q(H)\}.\]
\end{theorem}
\proof
In Section~\ref{subsec:quantum} we saw that
\[\chi_q(G \times H) \le \min\{\chi_q(G), \chi_q(H)\}.\]
Therefore we only need to show the reverse inequality. Suppose that $G$ and $H$ satisfy the conditions above. Recall from Lemma~\ref{lem:chain} that $\bartheta(K) \le \chi_q(K)$ for any graph $K$. Thus
\[\chi_q(G \times H) \ge \bartheta(G \times H) = \min\{\bartheta(G), \bartheta(H)\} = \min\{\chi_q(G), \chi_q(H)\},\]
by Theorem~\ref{thm:hedbartheta}.\qed

Recall that $\Omega_n$ is the graph with vertex set $\{\pm 1\}^n$ such that orthogonal vectors are adjacent. In Section~\ref{subsec:quantum}, we saw that these graphs exhibit exponential separation between $\chi_q$ and $\chi$ for $n$ a multiple of 4, and they have been central to the investigation of quantum chromatic number since its beginnings.

For $n$ odd, $\Omega_n$ is empty and thus $\chi_q(\Omega_n) = 1 = \bartheta(\Omega_n)$. For $n \equiv 2 \mod 4$, $\Omega_n$ is nonempty and bipartite, and thus $\chi_q(\Omega_n) = 2 = \bartheta(\Omega_n)$. For $n$ a multiple of 4, combining results from~\cite{Avis} and~\cite{qhomos} shows that $\chi_q(\Omega_n) = n = \bartheta(\Omega_n)$. Therefore, $\chi_q(\Omega_n) = \bartheta(\Omega_n)$ for all $n$ and thus we have the following corollary.

\begin{corollary}
For any $m,n \in \mathbb{N}$,
\[\chi_q(\Omega_m \times \Omega_n) = \min\{\chi_q(\Omega_m), \chi_q(\Omega_n)\}.\]
\end{corollary}

\section{Concluding Remarks}

We have shown that the $\chi_{vec}$, $\bartheta$, and $\chi_q$ versions of Sabidussi's theorem hold. We have also shown that the $\bartheta$ version of Hedetniemi's conjecture holds, the $\chi_{vec}$ version holds for 1-homogeneous graphs, and the $\chi_q$ version holds for graphs with strict vector chromatic number equal to quantum chromatic number. It is not surprising that we were more succesful with the analogs of Sabidussi's theorem, as this seems to be the easier of the two problems in general. However, we conjecture that the $\chi_{vec}$ and $\chi_q$ versions of Hedetniemi's conjecture hold in general.

With the similarity between $\chi_{vec}$ and $\bartheta$, it is worthwhile considering why the proof of Theorem~\ref{thm:hedbartheta} cannot be used to prove a version of Hedetniemi's conjecture for $\chi_{vec}$. The proof of Theorem~\ref{thm:hedbartheta} relies on the following three properties of $\bartheta$:
\begin{enumerate}
\item $\bartheta(G \cart H) = \max\{\bartheta(G), \bartheta(H)\}$
\item $\bartheta(G \boxtimes H) \ge \bartheta(G)\bartheta(H)$
\item $\bartheta(G \cup H) \le \bartheta(G)\bartheta(H)$
\end{enumerate}
Combining the last two gives that
\[\bartheta(G)\bartheta(H) \le \bartheta(G \boxtimes H) \le \bartheta(G \cart H) \bartheta(G \times H),\]
which along with the first proves the theorem. We noted after Theorem~\ref{thm:sabibartheta} that (1) also holds for $\chi_{vec}$, and it can be shown (using essentially the same proof as for $\bartheta$) that (2) holds for $\chi_{vec}$ as well. However, (3) is false for $\chi_{vec}$, as already shown by Schrijver in~\cite{Schrijver} (his $\theta'$ is equal to $\chi_{vec}$ of the complement). Of course this does not mean that a version of Hedetniemi's conjecture for $\chi_{vec}$ cannot be proved, but a different approach is needed.

We can consider the same analysis for $\chi_q$. Theorem~\ref{thm:qsabi} shows that (1) holds for $\chi_q$. Item (4) of Lemma~\ref{lem:qhomocart} concerning the disjunctive product shows that $\chi_q(G * H) \le \chi_q(G)\chi_q(H)$, and then the same trick used to prove Corollary~\ref{cor:barthetaunion} shows that (3) holds for $\chi_q$. This leaves (2), but it is not hard to see $\chi(C_5 \boxtimes C_5) = 5$ and thus
\[\chi_q(C_5 \boxtimes C_5) \le \chi(C_5 \boxtimes C_5) = 5 < 9 = \chi_q(C_5)^2.\]
Note that $\chi_q(C_5) = 3$ follows from the fact that $\chi_q(G) = 2$ if and only if $\chi(G) = 2$, which was proven in~\cite{qchrom}.

Of $\chi_{vec}$ and $\chi_q$, it seems that proving the analog of Hedetniemi's conjecture for the former should be more tractable. This is because one can use strong duality when working with $\chi_{vec}$, whereas $\chi_q$ is not known to have this property. On the other hand, finding a counterexample to the conjecture (if one exists) is also likely easier for $\chi_{vec}$ since it can be computed efficiently, and $\chi_q$ is not even known to be computable.

\renewcommand{\bibname}{References}

\bibliographystyle{plainurl}


 \end{document}